\documentclass[a4paper]{amsart}
\usepackage[english]{babel}
\newtheorem{theorem}{Theorem}[section]         
\newtheorem{lemma}[theorem]{Lemma}             
     
\newtheorem{proposition}[theorem]{Proposition} 
         
\newtheorem{remark}[theorem]{Remark}           
\newtheorem{definition}[theorem]{Definition}

\numberwithin{equation}{section}
\newcommand{\nc}{\newcommand}
\nc{\R}{\ensuremath{\mathbb{R}}}
\nc{\C}{\ensuremath{\mathbb{C}}}
\nc{\K}{\ensuremath{\mathbb{K}}}
\nc{\Z}{\ensuremath{\mathbb{Z}}}
\nc{\N}{\ensuremath{\mathbb{N}}}
\nc{\Mo}{\ensuremath{\mathbb{M}}}
\nc{\A}{\ensuremath{\mathcal {A}}}
\nc{\E}{\ensuremath{\mathcal {E}}}
\nc{\B}{\ensuremath{\mathfrak {B}}}
\nc{\Kg}{\ensuremath{\mathfrak {K}}}
\nc{\OS}{\operatorname{OS}}
\nc{\re}{\ensuremath{\mathcal R}}
\nc{\M}{\ensuremath{\mathcal M}}
\nc{\V}{\ensuremath{\mathcal V}}
\nc{\namelistlabel}[1]{\mbox{#1}\hfil}
\newenvironment{namelist}[1]{%
\begin{list}{}
    {
      
      \settowidth{\labelwidth}{#1}
      \setlength{\leftmargin}{1.1\labelwidth}
    }
  }{%
\end{list}}
\begin{document}
\title[Orlik\---Solomon algebras]
{A note on Tutte polynomials and Orlik\---Solomon algebras}
\thanks{2000 \emph{Mathematics Subject Classification}: 
05B35,  14F40, 32S22.\\
 \emph{Keywords and phrases}: arrangement of
 hyperplanes, 
 matroid, 
Orlik-Solomon algebra, Tutte polynomial.
}
\author{Raul Cordovil and David Forge}
\address{\newline 
Departamento de Matem\'atica,\newline 
Instituto Superior T\' ecnico \newline
 Av.~Rovisco Pais
 - 1049-001 Lisboa  - Portugal}
\email{cordovil@math.ist.utl.pt}
\thanks{The  first  author's research was 
supported in part by FCT (Portugal) through program POCTI and
 the project SAPIENS/36563/99.
The second  author's research  was  supported  by FCT trough the project 
SAPIENS/36563/99.}
\address{{}\newline
Laboratoire de Recherche en Informatique\newline 
Batiment 490
Universite Paris Sud\newline 
91405 Orsay Cedex -
France
}
\email{forge@lri.fr}
\begin{abstract}
 {Let $\A_{\C}=\{H_{1},\ldots, H_{n}\}$ be a (central) arrangement of 
hyperplanes in $\C^d$ and $\M(\A_{\C})$ the dependence matroid
of the linear forms $\{\theta_{H_i}\in (\C^d)^{*}:  
\mathrm{Ker}(\theta_{{H}_{i}})={H_i}\}$.
The Orlik-Solomon 
algebra
 $\OS(\M)$ of a matroid $\M$ is the   exterior algebra on the points 
modulo the ideal generated by
circuit boundaries. The graded algebra
$\OS(\M(\A_{\C}))$ is isomorphic to the 
cohomology
algebra  of the manifold\,
$\mathfrak{M}=\C^d\setminus\bigcup_{H\in
 \mathcal{A}_{\C}} H$. 
The Tutte polynomial $\operatorname{T}_\M(x,y)$ is a
powerful invariant  of the matroid $\M$. 
When $\M(\A_{\C})$ is a rank three  matroid and the
  $\theta_{H_i}$ are  complexifications of
real  linear forms, we will prove
  that  $\OS(\M)$
determines $\operatorname{T}_\M(x,y)$. This result solves partially a
conjecture of M. Falk.}
\end{abstract}
\thanks{Typeset by \AmS -\LaTeX}
\maketitle
\today
\section{Introduction}
Let $V$ be a vector space of dimension $d$  over some field $\K$.  
A (central) arrangement of 
hyperplanes in $V,$
$\mathcal{A}_{\mathbb {K}}=\{H_{1},\ldots, H_{n}\},$ is a
   finite listed set of codi\-mension one vector
subspaces. 
Given an arrangement $\A_{\K}$\, we suppose always 
selected a 
family of  linear forms 
$\{\theta_{{H_i}}:\, H_i\in \A_{\K},\, \theta_{H_i}\in V^{*},\,
\mathrm{Ker}(\theta_{{H}_{i}})={H_i}\}$. Let 
$\mathcal{M}(\mathcal{A}_{\mathbb {K}})$ be
the {\em dependence matroid}\/ determined by the vectors 
$\V :=\{\theta_{{H_i}}:\,\, H_i\in \A_{\K}\}:$
i.e.,   the matroid on the ground set $[n]$
which has for its independent sets the indices of  
linearly independent subsets of  $\V$.
A 
 matroid $\M$ is said to be \emph{realizable} over 
$\mathbb {K}$ if there is an arrangement
$\mathcal{A}_{\mathbb {K}}$ such that
$\M=\M(\mathcal{A}_{\mathbb {K}})$ 
 The manifold \,$\mathfrak{M}(\mathcal{A}_{\mathbb{C}})=
\C^d\setminus\bigcup_{H\in
 \mathcal{A}_{\mathbb{C}}} H$\, plays an important
role in the 
 Aomoto-Gelfand multivariable theory 
 of  hypergeometric functions (see \cite{OT2} for a recent 
introduction from the point of  view of
arrangement theory).  In  \cite{OT},  the
determination of the cohomology algebra
 $H^{*}\big(\mathfrak{M}(\mathcal{A}_{\C}); \C\big)$\, from the matroid
 $\mathcal{M}(\mathcal{A}_{\mathbb {C}}),$ is accomplished by first 
 defining the Orlik\---Solomon algebra $\OS(\A_{\C})$
 in terms of generators and 
 relators which depends only on the matroid 
$\mathcal{M}(\mathcal{A}_{\mathbb{C}}),$  and then by
showing that this graded commutative
 algebra is isomorphic to
$H^{*}\big(\mathfrak{M}(\mathcal{A}_{\mathbb{C}}); \C\big)$.
The Orlik-Solomon 
algebras  have been then
intensively studied. 
One of the central  questions is to find combinatorial 
invariants of $\OS$. Descriptions of the
developments from the late 1980s to the end of 1999,  
together with the contributions of many authors,
can be found in \cite{Falk1} or \cite{Yuz}.

Throughout this note  $\mathcal{M}$  denotes
 a (\emph{simple})  rank $r$ matroid without loops or 
parallel elements  on the ground set
$[n]$.  
Let $\mathfrak{C}({\mathcal{M}})$   be the set of 
\emph{circuits} (i.e.,  minimal dependent sets
with respect to inclusion)  of $\mathcal{M}$. 
 Let
 $\B(\M)$ be the set of \emph{bases} of the matroid.  
 The
\emph{Tutte polynomial}  of a matroid
$\M,$  denoted $T_\M(x,y),$ is the unique two-variable polynomial over
$\Z,$
 satisfying the following recursive 
properties for every $\ell\in [n],$ see \cite{Brylawski-Oxley}:
\begin{namelist}{xxx}
\item[$\,\,\circ$] $\operatorname{T}_\M(x,y)=\operatorname{T}_{\M\setminus
\ell}(x,y)+\operatorname{T}_{\M/\ell}(x,y)$\,  
if $\ell$ is neither a loop nor a coloop (isthmus),
\item[$\,\,\circ$] $\operatorname{T}_\M(x,y)=x
\operatorname{T}_{\M\setminus
\ell}(x,y)$
\,  if \, $\ell$\, is a coloop, 
\item[$\,\,\circ$] $\operatorname{T}_\M(x,y)=y
\operatorname{T}_{\M\setminus
\ell}(x,y)$\, 
if  \,$\ell$\, is a loop,
\item[$\,\,\circ$] $\operatorname{T}_\M(x,y)=1$ if the ground set of $\M$ is the emptyset.
\end{namelist}
In this note it is important to recall that   the Tutte polynomial of  a
rank three matroid $\M$ merely records the number of rank two flats  of
each cardinality.  More generally the Tutte polynomial of 
$\M$ can be reconstructed from
the multiset of isomorphism types of its hyperplanes, see 
\cite{Brylawski-Oxley}.
We refer to 
\cite{W1, W2}  
 as standard sources  for  matroids. \par
The authors wishes to thank a referee for his helpful suggestions.
\section{The main result}
Let  $E:=\{e_1, e_2,\dotsc,e_n\}$ be a fixed basis of 
the vector space $\C^n$. Consider the
\emph{graded exterior algebra} of the vector space\, 
$\bigoplus_{i=1}^n{\C}e_i:$
$$
\E:=\bigwedge(\bigoplus_{i=1}^{n}{\C}e_i)=
\E_0( = \C)\oplus\dotsb \oplus\E_\ell(=\sum\alpha_i e_{i_1}
 \dotsb\, e_{i_\ell},\, \alpha_{i}\in
\C)\,
\oplus\dotsb \oplus\E_n.
$$ 
When necessary we see the subset 
$\{i_1,\dotsc,i_\ell\}\subset [n]$ as  
 the ordered set $i_1<\dotsb<i_\ell$. Set
$e_X:=e_{i_1}\wedge e_{i_2}\wedge \dotsm\wedge e_{i_\ell},$ in particular 
set 
$e_{ij}:=e_i\wedge e_j$. By
convention
 set  $e_{\emptyset}:=1$. Let $\partial:\E\to \E$ be the unique
 morphism such that
$\partial(e_i\wedge e_j)=\partial(e_i)\wedge e_j-e_i\wedge\partial(e_j)$ 
and $\partial(e_i)=1,$
for all
$e_i, e_j \in E$ and $n \in \N$. 
 Let
$\Im(\mathcal{M})=\big\langle\partial(e_{\scriptstyle C}): C\in
\mathfrak {C}(\mathcal{M})\big\rangle$ be the  two-sided 
ideal of $\E$ generated by $\{\partial(e_{\scriptstyle C}): C\in
\mathfrak {C}(\mathcal{M})\}$.
\begin{definition}\label{def2}
{\em
The quotient $\E\big/ \Im(\M)$ is  termed the {\em\/ Orlik-So\-lo\-mon\/
$\C$-algebra} 
   of the matroid
 $\mathcal{M},$ denoted $\OS(\M)$. 
}
\end{definition}
The residue class in 
  $\OS(\M)$ 
  determined by the element $e_X\in \E$  is denoted  by $e_X$ as well.  
The algebra $\OS(\M)$ has the natural graduation 
$$\OS(\M)=\C\oplus\OS_1(=\E_1)\oplus\dotsm 
\oplus\OS_\ell(=\E_\ell/\Im(\M)\cap
\E_\ell)\oplus\dotsm \oplus \OS_r.$$
 We start by recalling some definitions and results from Falk 
\cite{Falk3, Falk1} about 
the resonance varieties of the $\C$-algebra $\OS(\M)$. 
We identify $\OS_1$ and $\C^n$.
Fix an element $e_\lambda=\sum_{i=1}^n \lambda_i e_i$ in $\OS _1$. 
Then left multiplication by 
$e_\lambda$ defines a map $\OS _p\rightarrow \OS_{p+1},$ 
which squares to zero. Thus it 
defines a cochain complex whose cohomology determines a stratification 
of the parameter space $\C^n$.  The $p^{th}$ resonance variety of 
$\OS$ is defined by 
\begin{equation*}
\re_p(\OS):=\{ \lambda\in \C^n\, |\, H^p(\OS,e_\lambda)\not=0\}.
\end{equation*}

\vspace*{2mm}
\noindent The $\re_p(\OS)$ are invariants of $\OS$.  They are the 
unions  of 
 subspaces of $\C^{n},$ see Corollary 3.8 in \cite{CO}.
 More is known about  $\re_1(\OS))$: it 
is a union of two kinds of linear subspaces of 
$\C^n$ with trivial intersection. The first are 
the \emph{local} components which correspond to flats of rank 2 of
 $k, k\geq 3,$
elements  and are of dimension $k-1$. The second are the \emph{non local}  
components which correspond to  some particular colorations 
of the elements  of  submatroids  of $\M,$ 
see \cite{Falk3}  for details.
A very interesting application of the 
variety $\re_{1}(\OS(\M))$ has been shown by M. Falk:
\begin{theorem}[Theorem 3.16 in \cite{Falk3}]\label{falk}
 Let $\M$ be a rank three matroid. 
Suppose every non local component of  $\re_{1}(\OS(\M))$ has dimension 
two.
Then   $\OS(\M)$ determines the Tutte polynomial
of $\M$.\qed
\end{theorem}
In relation to this theorem
Falk conjectured that if $\M$ has rank three then $\OS(\M)$ 
determines the Tutte polynomial
of $\M$. We prove this conjecture for matroids realizable over $\R$.
The proof is obtained simply by combining Theorem~\ref{falk} with 
a recent result of 
 Libgober and Yuzvinsky
and an old coloration result of Herzog and Kelly.
The next theorem of  Libgober and Yuzvinsky, which is 
fundamental in this note, makes more precise previous 
results of  Falk, see \cite{Falk3}.
 A \emph{regular} coloration of $\M([n])$ is a partition 
$\Pi$ of $[n]$ such that 
$\#\{\pi \in \Pi| \pi\cap X\not= \emptyset\}$ is equal to 1 or $|\Pi|$
for  all  rank two flats $X$ of $\M$.
\begin{theorem}[Theorem 3.9 and Remark 3.10 in \cite{LY}]\label{LY2}
Let $\Pi$ be a $k$-coloration,\\ $k\geq 2,$ of $\M$ which 
determines a non local component of\,
$\re_1(\OS(\M))$. Then $k\geq 3,$ $\Pi$ is regular 
and the dimension of the non local component
is $k-1$.\qed
\end{theorem}
We recall the following ``Sylvester type" result of Herzog and Kelly:
\begin{lemma}[Theorem 4.1 of \cite{HK}]\label{roubado}
Let $\M$ be a rank three
 matroid realizable over $\R$. Then it has no regular 
coloration with four or more colors.
\end{lemma}
\begin{proof}
Let $\A$ be an arrangement of lines in $\mathbb R^2$.
Let suppose that there exists a  regular coloration 
$\Pi$ of the matroid $\M(\A)$ with $k, k\geq 4,$ 
colors. The lines [resp.\, points] in $\A$ corresponding 
to elements [resp.\, rank two flats] of $\M(\A)$. 
Then all
multi-colored points of intersection of
$\A$ are intersection of at least 
$k$ lines of different color. Let  $P_0$ be a multi-colored 
point and $L_0\in \A, P_{0}\not \in L_{0},$ a line such that 
$d(P_0,L_0),$ the euclidean distance between $P_0$ and $L_0$, is minimum.
Let $L_1, L_2,\dots, L_k$ be $k$ lines through $P_0$ of 
different colors. At least 3 
of them are not of the color of $L_0$ and at least 2 of these 3 are 
not separated  by line $\ell$ perpendicular from $P_0$ 
to the line $L_0$. Let w.l.o.g  $P_0,L_0,L_1,
L_2$ be as in Figure 1. Then we have 
$d(L_0\cap L_2,L_1)< d(P_0,L_0)$ a contradiction.
\end{proof}
\begin{center}
\begin{picture}(-80,120)(100,-10)
\thicklines
\put(10,10) {\line(0,1){80}}
\put(10,10) {\line(3,1){180}}
\put(10,10) {\line(1,1){80}}
\put(10,10) {\line(0,-1){20}}
\put(10,10) {\line(-3,-1){20}}
\put(10,10) {\line(-1,-1){20}}
\put(50,50){\line(1,-3){7.9}}
\put(-10,50) {\line(1,0){200}}
\put(-5,12){$P_0$} \put(3,70){$\ell$}
\put(-25,50){$L_0$}\put(55,80){$L_2$}\put(160,75){$L_1$}
\put(10,10) {\circle*{8}}\put(50,50) {\circle*{8}}
\put(130,50){\circle*{8}}
\end{picture}
\put(20,-20)
{\textbf{Figure 1}}
\end{center}
\vspace{5mm}
Note that Lemma~\ref{roubado} is not true for matroids realizable over a
finite field. Indeed consider the affine plane  $\operatorname{AG}(2,
p^n),\, p^n\geq 4$. Choose a (maximal) set $L$ of $p^n$  ``parallel" 
(i.e., not intersecting two-by-two) lines of $\operatorname{AG}(2,
p^n)$. 
Color the points of every line of $L$ with the same color and two
different lines with two different colors. So, every line of the affine
plane not in $L$ is  $p^n$-colored  and the coloration is regular.
(See Example 4.5 in \cite{Falk3}, where a similar regular 3-coloration of 
the affine plane
$\operatorname{AG}(2,3)$ over $\Z_3$ is given.)
We do not know if the lemma remains true for line arrangements  in
$\C^2$.
Theorem \ref{questao} below solves
 partially   Conjecture~4.8 in \cite{Falk3}. It remains open the important
case of  matroids realizable over the complexes but not realizable over
the reals.
\begin{theorem}\label{questao}
If $\M$ is a matroid of rank at most three, realizable 
over $\R,$   then the Orlik-Solomon algebra $\OS(\M)$
determines the Tutte polynomial $\operatorname{T}_{\M}
(x,y)$.  
\end{theorem}
\begin{proof}
The Theorem is a consequence of Theorems\,\,\ref{falk}, 
\,\ref{LY2} and 
Lemma~\ref{roubado}.
\end{proof}
\section{Appendix}
From Lemma~\ref{roubado} and Theorem~\ref{LY2}
we know that in the rank 3 real case the non local components 
come only from regular 3-colorations, see Theorem 3.10 
in \cite{Falk3} for details. In \cite{EK,gru} we can find  the following 
regular 3-colorable real matroids $\M(\A)$ defined by line arrangements: 
\begin{namelist}{xxx}
\item[\,\,$\circ$]There is an infinite family of 
arrangements of $3k, k\geq 3,$ lines.
Take a regular $2k$-gon with its $k$ longest diagonals. 
Color the lines of the $2k$-gon 
alternatively  in two colors and the diagonals with a 
third color, see \cite {EK, gru}.
\item[\,\,$\circ$]
There is also the following simplicial arrangement, 
$\A_1(12),$ of 12 lines.Take the regular 6-gon together with its three 
longest and its three shortest diagonals.  Color 
the lines of the 6-gon cyclically in three colors
as well as the diagonals longest such that the vertices of the 6-gon
 are 3-colored. Then there is an unique coloration of the 
 three shortest diagonals such the 3-coloration is regular, 
 see
\cite{gru}.
\end{namelist}\par
Theorem~\ref{questao} does not hold for matroids of rank at least 4. 
Infinite families of counter-examples to Theorem~\ref{questao} 
appeared in \cite{Falk2}. 
We improve these  counter-examples to the class of ``connected matroids".
We make use of the Proposition~\ref{lemma} below. We say that 
the matroid $\M^{\prime}([n+1]),$ 
is a \emph{single-element} 
\emph{free extension} of its restriction  $\M^{\prime}([n+1])\mid [n]
:=\M([n]),$
if
$$\mathfrak{C}(\M([n+1]))=\mathfrak{C}(\M([n]))\cup
\{B\cup\{n+1\}:\, B\in \B(\M)\}.$$ We set $\M([n])\hookrightarrow 
\M^\prime([n+1])$ if $\M^\prime([n+1])$ is
free extension of $\M([n])$.
\begin{proposition}\label{lemma} 
Consider the two free extension $\M_1([n])\hookrightarrow 
\M_1^\prime([n+1])$ and
$\M_2([n])\hookrightarrow \M_2^\prime([n+1])$ and suppose 
that there is the graded isomorphism $\,\Phi:
\OS\big(\M_1([n])\big)\to 
\OS\big(\M_2([n])\big)$. Then there is a graded 
isomorphism
$$\,\Phi^\prime:
\OS\big(\M_1([n+1])\big)\to 
\OS\big(\M_2([n+1])\big)$$  extending $\,\Phi$ and such that  
$\Phi^\prime(e_{n+1})=e_{n+1}$. 
\end{proposition}
\begin{proof}  Set 
$\E^*:=\bigwedge (\bigoplus_{i=1}^{n+1}\C{e}_i)$.  We claim that: 
\begin{equation}\label{eq1}
\Im(\M_1')=\big(\Im(\M_1)\oplus\Im(\M_1)\wedge e_{n+1}\big)+
\big{\langle}
\partial(\E_{r+1}^{*})\big\rangle
\end{equation}
and
similarly for $\Im(\M'_2)$. Indeed
consider the ideal $\Delta_{r+1}$  of $\,\E^{n+1}$
$$ 
\Delta_{r+1}:=\Big{\langle}\Big\{\partial(e_{Y}): \,\, {Y}\in 
\binom{[n+1]}{r+1},\,\,\text{and}\,\,
(Y\setminus
\{n+1\})
\in   
\B(\M_1)\Big{\}}\Big{\rangle}.
$$
It is clear that $$\Im(\M_1)=
\big(\Im(\M_1)\oplus\Im(\M_1)\wedge
e_{n+1}\big)+\Delta_{r+1}.$$
Consider a subset $X^\prime\in  \binom{[n+1]}{r+1},$ and suppose that $
X^\prime\setminus
\{n+1\}=X^{\prime\prime}\not \in \B(\M_1)$.
Then  we have  $\partial(e_{X^{\prime}})\in
\Im(\M_1)\oplus\Im(\M_1)\wedge e_{n+1}$ and 
$\partial(e_{X^{\prime\prime}})\in
\Im(\M_1)$.  So
$$\big(\Im(\M_1\oplus\Im(\M_1)\wedge
e_{n+1}\big)+\big{\langle}
\partial(\E_{r+1}^{*})\big\rangle
\subset\big(\Im(\M_1)\oplus\Im(\M_1)\wedge
e_{n+1}\big)+\Delta_{r+1}.
$$
As the other inclusion is trivial  Equality $(\ref{eq1})$ follows.  The
restriction of the automorphism
$\Phi$ to $\mathrm{OS}_1(\M_1)(=\E_1)$  fixes a graded
automorphism
$\widetilde{\Phi}:\E\to \E,$ so
$\widetilde{\Phi}\big(\Im(\M_1)\big)=\Im(\M_2)$. Consider now the graded 
automorphism 
$\Xi:\E^{*}\to
\E^{*},$ such that  $\Xi\mid\E=\widetilde{\Phi}$ and\,
$\Xi(e_{n+1})=e_{n+1}$. It is clear that:
\begin{align*}
\Xi\big(\Im(\M_1')\big)=
\big(\Im(\M_2)\oplus\Im(\M_2)\wedge
e_{n+1}\big)+\big\langle
\Xi\big(\partial(\E^{*})\big)\big\rangle=\\
=\big(\Im(\M_2)\oplus\Im(\M_2)\wedge
e_{n+1}\big)+
\big\langle\partial\big(\Xi(\E^{*})\big)\big\rangle=
\Im(\M_2').
\end{align*}
\end{proof}
\begin{remark}
{\em If $\M([n])\hookrightarrow \M^\prime([n+1])$ the matroid 
$\M'([n+1])/n+1$  is called the
\emph{truncation} of
$\mathcal{M}([n])$. 
So we deduce from from Proposition~\ref{lemma} that isomorphism 
of Orlik\---Solomon algebras  
preserves  truncation, see Theorem 3.11  in \cite{Falk1}. 
Conversely we know that free extension by a point is
the same as adding a coloop and truncation. It is clear that isomorphism 
of
Orlik\---Solomon algebras is preserved
   under
directed sum with a coloop. 
So  Proposition~\ref{lemma} could be also deduced from Theorem 3.11  
in \cite{Falk1} and these two results are equivalent.
We leave the details to the reader.
}
\end{remark}
Consider now the two rank four matroids $\M_1$ and $\M_2,$  
given by the affine dependencies of the  seven
points in
$\R^3,$ as indicated in Figures 2 and 3 and characterized 
by the following properties:
\begin{namelist}{xxx}
\item[~$\circ$]The point 7 is in general position in both
 matroids.
\item[~$\circ$] Both  matroids have  exactly two lines 
with 3 points and these lines are the unique
lines with at least three elements.
\item[~$\circ$]  The
union of
the two lines with three elements has rank four in $\M_1$.
\item[~$\circ$]  The
union of
the two lines with three elements has rank three in $\M_2$.
\end{namelist}
The number of bases of the first matroid is 
$\operatorname{T}_{\M_3}(1,1)=27$
 and of the second matroid is    $\operatorname{T}_{\M_4}(1,1)=26$.   So 
we
have
$\operatorname{T}_{\M_1}(x,y)\not = \operatorname{T}_{\M_2}(x,y)$.
 We claim that
$\OS(\M_1) \cong \OS(\M_2)$. From Proposition~\ref{lemma} 
it is enough to prove that
$\OS(\M_1\setminus 7) \cong \OS(\M_2\setminus 7)$. 
Note that the matroids $\M_1\setminus 7$ and
$\M_2\setminus 7$ are  particular cases of the graphic matroids 
considered in \cite{Falk2}. (Make 
$\M_0$ and $C_n$ both  equal to triangles, in Figure 1 of 
\cite{Falk2}.) So the isomorphism follows from
 Theorem 1.1 of
\cite{Falk2}.\par
\vspace*{4cm}
\begin{picture}(0,0)(0,0)
\thicklines
\put(10,10) {\line(1,0){120}}
\put(10,10) {\line(3,1){60}}
\put(10,10) {\circle*{7}}\put(-1,10){3}
\put(70,30) {\circle*{7}}\put(68,17){4}
\put(70,30){\line(0,1){60}}
\put(130,10) {\circle*{7}}\put(135,10){6}
\put(130,10){\line(-3,1){60}}
\put(10,10) {\line(3,4){60}}
\put(40,50) {\circle*{7}}\put(29,50){2}
\put(70,90) {\circle*{7}}\put(59,90){1}
\put(70,90) {\line(3,-4){60}}
\put(100,20) {\circle*{7}}\put(107,20){5}
\put(60,60) {\circle*{7}}\put(59,46){7}
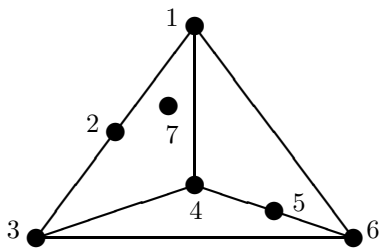
\put(-2,-15){\textbf{Figure 2}: The matroid $\M_1([7])$.}
\put(210,10) {\circle*{7}}\put(199,10){3}
\put(210,10) {\line(1,0){120}}
\put(210,10) {\line(3,1){60}}
\put(270,30){\line(0,1){60}}
\put(330,10){\line(-3,1){60}}
\put(210,10) {\line(3,4){60}}
\put(270,90) {\line(3,-4){60}}
\put(270,30) {\circle*{7}}\put(268,17){1}
\put(330,10) {\circle*{7}}\put(335,10){5}
\put(240,20) {\circle*{7}}\put(229,20){2}
\put(270,90) {\circle*{7}}\put(259,90){6}
\put(300,20) {\circle*{7}}\put(307,20){4}
\put(260,60) {\circle*{7}}\put(259,46){7}
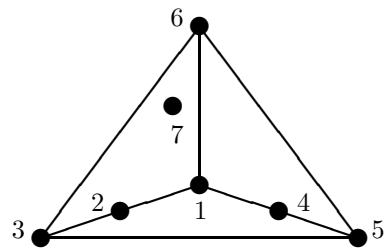
\put(199,-15){\textbf{Figure 3}: The matroid $\M_2([7])$.}
\end{picture}
\vspace{1.5cm}

\end{document}